\documentclass[11pt]{amsart}
\usepackage{amsmath,amssymb,amsthm}
\usepackage[latin1]{inputenc}
\usepackage{version,tabularx,multicol}
\usepackage{graphicx,float}

\newcommand{\reff}[1]{(\ref{#1})}

\theoremstyle{plain}
\newtheorem{theo}{Theorem}[section]

\newtheorem{prop}[theo]{Proposition}

\newtheorem{lem}[theo]{Lemma}

\theoremstyle{remark}
\newtheorem{rem}[theo]{Remark}

\newcommand{\ddelta}{\Delta}
\newcommand{\ca}{{\mathcal A}}

\newcommand{\cc}{{\mathcal C}}

\newcommand{\cf}{{\mathcal F}}
\newcommand{\cg}{{\mathcal G}}
\newcommand{\ch}{{\mathcal H}}

\newcommand{\cw}{{\mathcal W}}

\newcommand{\D}{{\mathbb D}}
\newcommand{\E}{{\mathbb E}}

\newcommand{\N}{{\mathbb N}}
\renewcommand{\P}{{\mathbb P}}

\newcommand{\R}{{\mathbb R}}

\newcommand{\W}{{\mathbb W}}

\newcommand{\rP}{{\rm P}}

\newcommand{\ind}{{\bf 1}}

\newcommand{\inv}[1]{\mathop{\frac{1}{ #1}}\nolimits}
\newcommand{\expp}[1]{\mathop {\mathrm{e}^{ #1}}}

\begin{document}

\title{Height process for super-critical continuous state branching
  process} 

\date{\today}

\author{Jean-François Delmas}

\address{CERMICS, \'Ecole Nationale des Ponts et Chaussées, ParisTech,
  6-8 av. Blaise Pascal, 
  Champs-sur-Marne, 77455 Marne La Vallée, France.}

\email{delmas@cermics.enpc.fr}

\thanks{This research  was partially supported by 
  NSERC Discovery Grants of the Probability group at Univ. of British Columbia.}

\begin{abstract}
  We  define  the height  process  for  super-critical continuous  state
  branching processes with quadratic  branching mechanism. It appears as
  a projective  limit of Brownian motions with  positive drift reflected
  at 0  and $a>0$ as $a$ goes  to infinity.  Then we  extend the pruning
  procedure of branching processes to the super-critical case. This give
  a  complete  duality picture  between  pruning  and size  proportional
  immigration for quadratic continuous state branching processes.
\end{abstract}



\maketitle

\section{Introduction}

Continuous
state branching  process (CB) appears as the limit of Galton-Watson
processes, see  \cite{f:dpg} for the quadratic branching mechanism and
\cite{l:lsbp} in the general case. We shall be interested in a CB,
$Z^\theta=(Z^\theta_r, r\geq 0)$,  with
quadratic branching mechanism, $\psi_\theta$, 
\begin{equation}
   \label{eq:def-psi}
   \psi_\theta(u)= 2u^2+4 \theta u,\quad u\geq 0,
\end{equation}
for  a  given  parameter  $\theta\in  \R$.   The process $Z^\theta$ is
a continuous  Markov process taking values in $\R_+$ 
such that for all $r\geq 0$, $\lambda\geq 0$, $x\geq 0$, 
\begin{equation}
   \label{eq:def-Z_q}
\E[\expp{-\lambda Z_r^\theta}|Z^\theta_0=x]=\expp{-x u^\theta(\lambda,r)},
\end{equation}
where $u^\theta$ is the only non-negative solution of the differential
equation 
\begin{equation}
   \label{eq:def-u_q}
 u'(r)+ \psi_\theta(u(r))=0, \quad r\geq 0,  \quad \text{and}\quad
 u(0)=\lambda. 
\end{equation}
(In fact the  general quadratic branching is of  the form $\displaystyle
\psi(u)=2\alpha  u^2   +  4\alpha   \theta  u$,  with   $\alpha>0$.  The
corresponding CB is distributed as $\displaystyle ( Z^\theta_{\alpha r},
r\geq 0)$.  Up  to this time scaling, we see it  is enough to consider
the case $\alpha=1$.)

The quantity $Z^\theta_r$ can be thought  as the ``size'' at time $r$ of
a   population  of   individual  with   infinitesimal  mass   and  whose
reproduction mechanism  is characterized by  $\psi_\theta$.  The process
$Z^\theta$ is  called critical (i.e.   constant in mean)  if $\theta=0$,
sub-critical (i.e.   with exponential decay for the  mean) if $\theta>0$
and  super-critical (i.e.   with  exponential growth  for  the mean)  if
$\theta<0$.   In the  critical or  sub-critical case,  one can  code the
genealogy associated  to $Z^\theta$ using the  so-called height process,
$H^\theta=(H^\theta_t,  t\geq 0)$,  see \cite{lg:betmbp}  for $\theta=0$
and \cite{dlg:rtlpsbp} in a more general setting.  The height process is
the  limit   of  the  contour   processes  associated  to   sequence  of
Galton-Watson   trees  which   converge   to  $Z^\theta$.    Intuitively
$H^\theta_t$  is  the  genealogy  of  individual with  label  $t$  in  a
continuous  branching  process.   The  ``size''  of  the  population  of
individuals  with   label  less  than   $t$  and  which  are   alive  at
``generation'' $r$ is given by the local time of $H^\theta$ at level $r$
up to time $t$: $L_r^\theta(t)$.  To consider an initial population with
``size''  $x>0$, it  is  enough to  look  at the  height  process up  to
$T^\theta=\inf\{t>0;           L^\theta_0(t)=x\}$.           Intuitively
$L^\theta_r(T^\theta)$  gives the  ``size'' of  the population at
generation $r$ coming from an initial 
population with  ``size'' $x$.  In particular,
one      expect     the      height     process   
$L^\theta=(L^\theta_r(T^\theta), r\geq  0)$ to be  a CB started  at $x$ with
branching mechanism $\psi_\theta$.

For  $\theta=0$,  $H^\theta$ is  distributed as  the absolute
value of a  Brownian motion and the second  Ray-Knight theorem gives the
process $L^0$  is indeed a CB  with branching  mechanism $\psi_0$
started from $x$.  For  $\theta>0$,  $H^\theta$ is  distributed as the
Brownian motion with drift $-2\theta$ and reflected above $0$. 
Using Girsanov  theorem or the more general framework
developed  in  \cite{dlg:rtlpsbp},  it   is  easy  to  check  that  for
$\theta\geq  0$,  the process  $L^\theta$ is indeed a  CB  with
branching  mechanism $\psi_\theta$  started from  $x$.

Our aim  is to extend those  results to the case  $\theta<0$.  One would
like  to  consider  a  reflected  Brownian motion  with  positive  drift
$-2\theta$.   However, in this  case $T^\theta$  might be  infinite.  In
fact we  have that  $\P(T^\theta<\infty )$ is  equal to  the probability
that a CB with branching  mechanism $\psi_\theta$ become extinct that is
$1$ if $\theta\geq 0$ or $\displaystyle \expp{2 x\theta}$ if $\theta<0$,
see  \cite{g:abctcsbp}.   Intuitively, if  $T^\theta=\infty  $ it  means
there are  individuals alive  at generation $r=\infty  $, and  the height
process $H^\theta$  describes only (part of)  the lineage of  one of the
individuals alive  at time  $\infty $.  To  circumvent this  problem, we
chose  to consider the  height process  associated not  to the  whole CB
process but only up to a generation $a$.  In a discrete setting, we would
consider a (super-critical)  Galton-Watson process and its corresponding
discrete tree,  and would  cut the tree  above a given genealogy,  and would
look at the discrete  height process of this finite tree.  Following
the procedure in \cite{dlg:rtlpsbp}, one would expect the height process
of the discrete  tree to converge to a  Brownian motion
(with drift)  reflected at 0 and $a$. 

This  intuition lead us  to consider  for $a>0$  a Brownian  motion with
drift   $-2\theta$   and   reflected    above   $0$   and   below   $a$,
$H^{\theta,a}=(H^{\theta,a}_t, t\geq  0)$.  In section  \ref{sec:hp}, we
first  check that the  family $(H^{\theta,a},  a>0)$ can  be built  in a
consistent way. Let $\cc^c$ be the set of continuous function defined on
$\R_+$ taking values  in $[0,c]$. We define the  projection from $\cc^a$
to $\cc^b$, $\pi_{a,b}$, by $\pi_{a,b}(\varphi)(t)=\varphi(C_\varphi(t))$
for $t\geq 0$, where  $\varphi\in \cc^a$ and $C_\varphi(t)= \inf\{ r\geq
0; \int_0^r \ind_{\{\varphi(s)\leq b\}}\; ds>t \}$ is the inverse of the
time spent by $\varphi$ below $b$.  By construction, we have for $a>b>c$
that  $\pi_{a,c}=\pi_{a,b}\circ \pi_{b,c}$.   As the  process $\pi_{a,b}
(H^{\theta,a})    $   is    distributed    as   $H^{\theta,b}$    (Lemma
\reff{lem:proj}), this compatibility relation implies the existence of a
projective   limit  $\ch^\theta=(\ch^{\theta,a},a\geq   0)$   such  that
$\ch^{\theta,a}$  is distributed  as  $H^{\theta,a}$ and  $\displaystyle
\pi_{a,b} (\ch^{\theta,a})=\ch^{\theta,b}$.   We shall call $\ch^\theta$
the height process of the quadratic branching process. It is defined for
$\theta\in  \R$.   We  can  consider  $Z^\theta_r$  the  local  time  of
$\ch^{\theta,a}$ at  level $r$  up to  the hitting time  of $x$  for the
local   time  of  $\ch^{\theta,a}$   at  level   $0$.  Because   of  the
compatibility relation,  we shall see that $Z^\theta_r$  does not depend
on  $a$, as  soon as  $a\geq  r$.  We  prove a  Ray-Knight theorem  for
$\ch^\theta$: $(Z^\theta_r,  r\geq 0)$ is a CB  with branching mechanism
$\psi_\theta$  (see Theorem  \ref{th:RK-th}). Of  course we  recover the
critical  cases   and  sub-critical   cases,  see  comments   of  Remark
\ref{rem:ch=cst}.  The   proof  relies  on  Girsanov   theorem  and  the
Ray-Knight theorem for $\theta=0$.

Following  \cite{lg:sbprspde},  we  can  add  a spatial  motion  to  the
individuals  to get a  super-critical Brownian  snake. Taking  a Poisson
process as a spatial motion,  this allows to adapt the pruning procedure
developed  in \cite{as:rbsd,as:psf}  (see also  \cite{ad:falp}  for more
general critical  or sub-critical branching mechanism)  for the critical
case  to the  super-critical case.   This  procedure gives  a nice  path
transformation   to  get   $\ch^{\theta}$   from  $\ch^{\theta'}$   when
$\theta>\theta'$  belong to  $\R$,  see Proposition  \ref{prop:pruning}.
Using this pruning transformation and the Ray-Knight theorem, we can get
$Z^\theta$ from  $Z^{\theta'}$ for any $\theta>\theta'$ (this result is
new for $\theta'<0$).   Notice that a
size     proportional    immigration    procedure,     introduced    in
\cite{ad:cbpimcsbpi} in  a more  general setting, allows  to reconstruct
$Z^{\theta'}$ from  $Z^\theta$. Our  result complete the  description of
the  duality  between  size  proportional immigration  and  pruning  for
quadratic branching mechanisms.

The paper is organized as follows.  In Section \ref{sec:hp}, we check the
compatibility relation in order to define the height process in the
super-critical case. Section \ref{sec:RK} is devoted to the proof of the
Ray-Knight theorem. In Section \ref{sec:BS}, following
\cite{lg:sbprspde, lg:bssd} we define the Brownian snake for
super-critical branching mechanism. The pruning procedure is developed
is Section \ref{sec:p}. 

\section{Height process for quadratic branching process} 
\label{sec:hp}
We assume that $\theta\in \R$.  Let $H^{\theta,a}=(H^{\theta,a}_t, t\geq
0)$ be a  Brownian motion  with drift $-2\theta$  reflected in  $[0,a]$ and
started  at $0$.  This  process can  be constructed  using a  version of
Skorohod's equation (see \cite{t:sderbccr}). This is the unique solution
of the stochastic differential equation:
\begin{equation}
   \label{eq:sde}
dY_t=d\beta_t -2\theta dt +\inv{2} d L_0(t) -\inv{2} dL_a (t), \quad
Y_0=0,
\end{equation}
where $L_y(t)$ is the local time of $Y$ at level $y$ up to time $t$ and
$(\beta_t, t\geq 0)$ is a standard Brownian motion. 

We first check  that the family $(H^{\theta,a}, a>0)$ can  be built in a
consistent way. Let $\cc^c$ be the set of continuous function defined on
$\R_+$ taking values in $[0,c]$. Let $a>b>0$. For $\varphi\in \cc^a$, we
consider the  time spent below level  $b$ up to  time $t$: $A_t=\int_0^t
\ind_{\{\varphi(s)\leq  b\}} \;  ds$  and its  right continuous  inverse
$C_\varphi(t)=\inf\{ r\geq  0; A_r>t\}$, with the  convention that $\inf
\emptyset=\infty $  and $\varphi(\infty  )=b$. We define  the projection
from        $\cc^a$       to       $\cc^b$,        $\pi_{a,b}$,       by
$\pi_{a,b}(\varphi)=\varphi\circ C_\varphi$.

\begin{lem}
\label{lem:proj}
Let $a>b>0$.    The process $\pi_{a,b} (H^{\theta,a}) $ is distributed as
   $H^{\theta,b}$. 
\end{lem}

\begin{proof}
  For  convenience  we shall  write  $H^{}$  instead of  $H^{\theta,a}$.
  Notice that $H$ solves \reff{eq:sde}.  Let  $(\cg_t,  t\geq  0)$  be  the
  filtration generated by the Brownian motion $\beta$, completed the usual way.
Let $\displaystyle  A_t= \int_0^t  \ind_{\{H^{}_s\leq b\}}  \; ds$
and $\displaystyle C(t)=\inf\{ r\geq  0;
A_r >t\}$.  Notice that a.s.
the stopping time  $C(t)$ is finite.  Let $L_r^{}(t)$  be the local time
of $H^{}$  at level $r$  up to time  $t$.  We set  $\displaystyle \tilde
H_t= H^{}_{C(t)}$. Using \reff{eq:sde}, we get
\begin{align*}
   \tilde H_t
&=\int_0^{C(t)} (d \beta_s - 2 \theta ds) + \inv{2} L_0^{}(C(t)) -\inv{2}
L_a^{}(C(t)) \\
&=\int_0^{C(t)} \ind_{\{H^{}_s\leq  b\}} (d \beta_s - 2 \theta ds) + \inv{2} L_0^{}(C(t)) -\inv{2}
L_a^{}(C(t))  + \int_0^{C(t)} \ind_{\{H^{}_s>  b\}} (d \beta_s -
2 \theta ds).
\end{align*}
Since  $A$  is continuous,  ,  by  construction  we have  $\displaystyle
\int_0^{C(t)} \ind_{\{H^{}_s\leq b\}} ds=\int_0^{C(t)} dA_s=A_{C(t)}=t$.
Notice  that $\beta'=(\beta'_t,t\geq 0)$,  where $\beta'_t=\displaystyle
\int_0^{C(t)}  \ind_{\{H^{}_s\leq  b\}}  d  \beta_s$,  is  a  continuous
martingale (with respect to the filtration $(\cf_{C(t)}, t\geq 0)$). Its
bracket   is   given  by   $\langle   \beta'  \rangle_t   =\displaystyle
\int_0^{C(t)} \ind_{\{H^{}_s\leq b\}} ds=  t$.  Therefore, $\beta'$ is a
Brownian motion. On the other hand,  by Tanaka formula, we have a.s. for
$r\geq 0$,
\[
(H^{}_r - b)^+=
(H^{}_0 - b)^+ + \int_0^r \ind_{\{ H^{}_s> b\}}
  dH^{}_s+ \inv{2} L_b^{}(r),
\]
where $x^+=\max(x,0)$.  Since $H^{}_0=0$ and
$H^{}_{C(t)}\in [0,b]$, we get
with $r=C(t)$ that 
\[
\int_0^{C(t)} \ind_{\{ H^{}_s> b\}}
  dH^{}_s+ \inv{2} L_b^{}(C(t))=0.
\] 
Use 
\[
 \int_0^{r} \ind_{\{ H^{}_s> b \}}
  dH^{}_s=\int_0^{r} \ind_{\{ H^{}_s> b \}}  (d \beta_s -
2 \theta ds) -\inv{2} L^{}_a(r)
\]
to get 
\[
-\inv{2} L^{}_a(C(t))+ \int_0^{C(t)} \ind_{\{H^{}_s>  b\}} (d \beta_s -
2 \theta ds)= -\inv{2} L_b^{}(C(t)).
\]
Therefore, we have 
\[
\tilde H_t=\beta'_t-2\theta t  +\inv{2}
L_0^{}(C(t)) -  \inv{2} L_b^{}(C(t)),
\]
where $\beta'$ is a Brownian motion. 
Notice the function $K$ defined for $t\in \R_+$ by $K(t)=\inv{2}
L_0^{}(C(t)) -  \inv{2} L_b^{}(C(t))$ is continuous with bounded
variation such that $K(0)=0$. Furthermore we have 
\[
\int_0^\infty \ind_{\{\tilde H_s \not\in \{0,b\}\}} d|K|(t)=0
\]
and 
\[
dK(t)= \ind_{\{\tilde H_t=0\}} d|K|(t) - \ind_{\{\tilde H_t=b\}} d|K|(t).
\]
Since $\tilde H$ is a continuous function taking values in $[0,b]$, we
deduce from theorem 2.1 in \cite{t:sderbccr}, that $\tilde H$ is a Brownian motion with drift $-2\theta$ reflected in $[0,b]$ started
at $0$. Henceforth, it is distributed as $H^{\theta,b}$.  

\end{proof}

Notice    that   by    construction   we    have   for    $a>b>c$   that
$\pi_{a,c}=\pi_{a,b}\circ  \pi_{b,c}$. Let  $\mu_r$ denotes  the  law of
$H^{\theta,r}$ for $r\geq 0$.   Lemma \ref{lem:proj} entails that $\mu_a
\circ (\pi_{a,b})^{-1}=\mu_b$.  This compatibility relation  implies the
existence  of a  projective limit  $\ch^\theta=(\ch^{\theta,a},a\geq 0)$
such that $\ch^{\theta,a}$ is distributed as $H^{\theta,a}$ and
\begin{equation}
   \label{eq:PiH}
\pi_{a,b} (\ch^{\theta,a})=\ch^{\theta,b}. 
\end{equation}
We will call $\ch^\theta$ the height process of the quadratic branching
process. 

\begin{rem}
\label{rem:coincider}
  If there exists $a\geq b>  0$ s.t. $\ch^{\theta,a}$ does not reach $b$
  on  $[0,t]$,  then we  have  that  a.s.  $\ch^{\theta,c}$ coincide  on
  $[0,t]$ for all $c\geq b$. 
\end{rem}

\section{Ray-Knight theorem for reflected Brownian motion with drift}
\label{sec:RK}

Let $L^{\theta,a}_r(t)$ be the local time of $\ch^{\theta,a}$ at level
$r$ up to time $t$. For $x>0$ we define  
\begin{equation}
   \label{eq:def-T}
T^{\theta,a}_x=\inf \{ t\geq 0; L_0^{\theta,a}(t)>x\},
\end{equation}
with the convention $\inf \emptyset=\infty $. Let $r\geq 0$. Notice that
equation \reff{eq:PiH} implies that for all $a\geq b \geq r$,
\[
L^{\theta,a}_r(T^{\theta,a}_x)=L^{\theta,b}_r(T^{\theta,b}_x). 
\]
We shall denote this common value by $Z^\theta_r$. We write
$\ch^{\theta,a,(x)} =(\ch^{\theta,a}_t, t\in [0, T^{\theta,a}_x])$ and
we call 
$\ch^{\theta,(x)}=(\ch^{\theta,a,(x)}  ,
a\geq 0)$ the height process associated to $Z^\theta=(Z^\theta_r, r\geq
0)$. 

We can now formulate
the Ray-Knight theorem.

\begin{theo}
   \label{th:RK-th}
The process $Z^\theta$ is a CB with branching
mechanism $\psi_\theta$. 
\end{theo}

\begin{rem}
\label{rem:ch=cst}
On the event  that $Z^\theta $ become extinct, there  exists a level $r$
s.t.  $Z^\theta_r=0$,  that  is  $L^{\theta,b}_r(T^{\theta,b}_x)=0$  for
$b\geq   r$.    From   Remark   \ref{rem:coincider},  we   deduce   that
$(\ch^{\theta,a}_t,t\in [0,T^{\theta,a}_x]) $  does not depend on $a>r$.
In  the  sub-critical  or  critical  case  (i.e.  $\theta\geq  0$),  the
extinction is almost sure. Thus the process $(\ch^{\theta,a}_t, t\in [0,
T^{\theta,a}_x])$ is constant for $a$ large enough. It is distributed as
a Brownian motion with drift $-2\theta$ reflected above $0$ stopped when
its  local  time  at  level  $0$  reaches $x$.  In  this  case,  Theorem
\ref{th:RK-th}   correspond  to  the   usual  Ray-Knight   theorem  (see
\cite{ry:mb},  chap.  XI.2   for  $\theta=0$  and  \cite{w:bprtsbm}  for
$\theta>0$).
\end{rem}

\begin{proof}
  Let $a>0$ and  $x>0$ be fixed. To be concise,  we write $H^\theta$ for
  $H^{\theta,a}$  and  $T^\theta=\inf  \{t\geq 0;  L_0^{\theta}(t)>x\}$,
  where $L^\theta_r(t)$ is the local  time of $H^\theta$ at level $r$ up
  to time $t$.  Notice $T^\theta$ is finite a.s. Let $g$ be a continuous
  function taking values in $\R_+$.  By monotone convergence, we have
\[
\E\left[\expp{    -\int_0^{T^\theta}     g(H^\theta_s)    ds}    \right]
=\lim_{n\rightarrow  \infty } \E\left[\expp{  -\int_0^{T^\theta\wedge n}
    g(H^\theta_s) ds} \right].
\]
Using Girsanov theorem and the fact that $H^\theta$ solves
\reff{eq:sde}, where $L_0$ and $L_a$ are continuous adapted functionals of
$(\beta_t -2\theta t, t\geq 0)$ (this is a consequence of Theorem 2.1
in \cite{t:sderbccr}), we get that 
\[
\E\left[\expp{  -\int_0^{T^\theta\wedge n}
    g(H^\theta_s) ds} \right]
=\E\left[\expp{ -2\theta \beta_{T^0\wedge n}  - 2 \theta^2 (T^0\wedge
    n) } \expp{  -\int_0^{T^0\wedge n}
    g(H^0_s) ds} \right].
\]
Since $H^0$ solves equation \reff{eq:sde} with $\theta=0$, we deduce
that 
\begin{equation}
   \label{eq:beta_T}
\beta_{T^0\wedge n} = H^0_{T^0\wedge n} - \inv{2} L_0^0(T^0\wedge n)
+ \inv{2}  L_a^0(T^0\wedge n).
\end{equation} 
Since $L_0^0(T^0)=x$,  we have $\displaystyle 
\beta_{T^0\wedge n} \geq - \inv{2} L_0^0(T^0\wedge n)\geq - \inv{2}
L_0^0(T^0)=-\frac{x}{2}$.
By monotone convergence, we get 
\[
\lim_{n\rightarrow  \infty } \E\left[\expp{ -2\theta \beta_{T^0\wedge n}
    - 2 \theta^2 (T^0\wedge 
    n) } \expp{  -\int_0^{T^0\wedge n}
    g(H^0_s) ds} \right]
=\E\left[\expp{ -2\theta \beta_{T^0}  - 2 \theta^2 T^0 } \expp{
    -\int_0^{T^0}     g(H^0_s) ds} \right]. 
\]
We write $Z_r=L^0_r(T^0)$ for $r\in [0,a]$. 
Notice that \reff{eq:beta_T} implies 
\[
\beta_{T^0}=  \inv{2} Z_a - \frac{x}{2}. 
\]
This and the occupation time formula for $H^0$ implies that $T^0=\int_0^a Z_r \;
dr$ and 
\[
\E\left[\expp{ -2\theta \beta_{T^0}  - 2 \theta^2 T^0 } \expp{
    -\int_0^{T^0}     g(H^0_s) ds} \right]
=\E\left[\expp{ \theta x - \theta  Z_a  - 2 \theta^2 \int_0^a Z_r \;
dr } \expp{
    -\int_0^{a}     g(r)Z_r \;  dr} \right].
\]
This leads to 
\begin{equation}
   \label{eq:GirH}
\E\left[\expp{    -\int_0^{T^\theta}     g(H^\theta_s)    ds}    \right]
= \E\left[\expp{ \theta x - \theta  Z_a  - 2 \theta^2 \int_0^a Z_r \;
dr } \expp{
    -\int_0^{a}     g(r)Z_r \;  dr} \right].
\end{equation}
Use the time occupation formula for $\ch^{\theta,a,(x)}$ (which
is distributed as $(H^\theta_s,s\in [0,T^\theta])$) to  get 
\begin{equation}
   \label{eq:Zq}
\E\left[ \expp{
    -\int_0^{a}     g(r) Z_r^\theta \;  dr} \right]
= \E\left[\expp{ \theta x - \theta  Z_a  - 2 \theta^2 \int_0^a Z_r \;
dr } \expp{
    -\int_0^{a}     g(r)Z_r \;  dr} \right].
\end{equation}

The Ray-Knight theorem implies  that $Z=(Z_r=L^0_r(T^0), r\in [0,a])$ is
distributed as the square of 0-dimensional Bessel process started at $x$
up to time $a$. In particular it is the unique strong solution of
\[
d\hat Y_t= 2 \sqrt{\hat Y_t}\; dW_t \quad t\in [0,a], \quad \hat Y_0=x,
\]
where $(W_t, t\geq 0)$ is a standard Brownian motion in $\R$. We deduce 
\[
\E\left[\expp{ \theta x - \theta  Z_a  - 2 \theta^2 \int_0^a Z_r \;
dr } \expp{
    -\int_0^{a}     g(r)Z_r \;  dr} \right]
=\E\left[\expp{- 2 \theta \int_0^a  \sqrt{\hat Y_t} dW_t   - 
2\theta^2 \int_0^a \hat Y_r \;
dr } \expp{
    -\int_0^{a}     g(r)\hat Y_r \;  dr} \right].
\]
Notice that $M=(M_t,t\geq 0)$, where $\displaystyle  M_t=\expp{- 2
  \theta \int_0^t
  \sqrt{\hat Y_r} dW_r   -  
2\theta^2 \int_0^t \hat Y_r \;
dr }$, define a local  martingale. It is in fact a martingale (see
section 6 in \cite{py:dbb}). 
Using Girsanov theorem again, we get that 
\[
\E\left[\expp{- 2 \theta \int_0^a  \sqrt{\hat Y_t} dW_t   - 
2\theta^2 \int_0^a \hat Y_r \;
dr } \expp{
    -\int_0^{a}     g(r)\hat Y_r \;  dr} \right]
=\E\left[ \expp{
    -\int_0^{a}     g(r)\hat Y_r^\theta \;  dr} \right],
\]
where $\hat Y^\theta$ is the
 unique strong solution of the stochastic differential equation 
\[
d\hat Y_t^\theta= 2 \sqrt{\hat Y_t}\; dW_t^\theta -4\theta \hat Y^\theta_t
\; dt \quad t\in [0,a], \quad \hat Y_0^\theta=x,
\]
where $W^\theta$ is a standard Brownian motion.
In conclusion we get 
\[
\E\left[ \expp{
    -\int_0^{a}     g(r) Z_r^\theta \;  dr} \right]
=\E\left[ \expp{
    -\int_0^{a}     g(r)\hat Y_r^\theta \;  dr} \right].
\]
We  deduce  that  $(Z_r^\theta,r\in  [0,a])$ is  distributed  as  $(\hat
Y_r^\theta ,r\in [0,a])$  for all $a>0$. In particular,  $Z^\theta$ is a
continuous Markov process. 
Recall $u^\theta$  defined by \reff{eq:def-u_q}. 
We have $\displaystyle u^0(\lambda,s)=\lambda/(1+2\lambda s)$ and for
$\theta\neq 0$
\begin{equation}
   \label{eq:form-uq}
 u^\theta (\lambda,s)=\frac{\lambda \expp{-4
  \theta t} }{1+\lambda (2\theta)^{-1} ( 1-\expp{-4\theta t})}. 
 \end{equation} 
From \reff{eq:Zq} we deduce that for
$\lambda\geq 0$ 
\[
\E\left[ \expp{
    -\lambda  Z_a^\theta } \right]
= \E\left[\expp{ \theta x - (\theta+\lambda)  Z_a  - 2 \theta^2 \int_0^a Z_r \;
dr } \right].
\]
Thanks to formula (2.k) in \cite{py:dbb}, we check the right hand side
is equal to $\displaystyle \expp{- x u^{\theta}(\lambda,a)}$. This
implies that $Z^\theta$ is a CB with branching mechanism $\psi_\theta$. 
\end{proof}

\section{The Brownian snake}
\label{sec:BS}
\subsection{Definition}

We refer  to \cite{dlg:rtlpsbp}, section 4.1.1, for  the construction of
the snake with a fixed lifetime process. Let $\xi$ a Markov process with
càdlàg paths and values in a Polish space $E$, whose topology is defined
by a metric $\delta$. We  assume $\xi$ has no fixed discontinuities. Let
$\rP_y$ denote the law of $\xi$ started at $y\in E$. The law of $\xi$ is
called the spatial motion. For $y\in  E$, let $\W_y$ be the space of all
$E$-valued killed paths started at $y$. An element of $\W_y$ is a càdlàg
mapping $w: [0,\zeta) \rightarrow  E$ s.t. $w(0)=y$. $\zeta\in (0,\infty
)$ is called the lifetime of $w$.
By convention the point $y$ is considered as the
path with zero lifetime, and is added to $\W_y$. 
The space
$\W=\cup_{y\in E} \W_y$, equipped with the distance defined in
\cite{dlg:rtlpsbp}, section 4.1.1, 
is  Polish. For $w\in
\cw$,  we define $\hat  w=w(\zeta-)$ if  the limit  exists and  $\hat w
=\Delta$ otherwise, where $\Delta$ is a cemetery point added to $E$.

Mimicking
the  proof  of  proposition  4.1.1  in  \cite{dlg:rtlpsbp},  for  $a>0$,
$\theta\in    \R$,    there    exists    a   càdlàg    Markov    process
$W^{\theta,a}=(W_s^{\theta,a}, s\geq 0)$ taking values in $ \W_y$ s.t.
\begin{itemize}
   \item If $\zeta_s^{\theta,a}$ denotes the lifetime of
     $W_s^{\theta,a}$, then $\zeta^{\theta,a}=(\zeta_s^{\theta,a}, s\geq
     0)$  is
     distributed as $H^{\theta,a}$. 
   \item   Let   $s\geq   0$.   Conditionally   on   $\zeta^{\theta,a}$,
     $W^{\theta,a}_s$    is    distributed     as    $\xi$ under
     $\rP_y$     on    $[0,
     \zeta^{\theta,a}_s)$. Notice that a.s. 
     $\hat W^{\theta,a}_s=W^{\theta,a}_s(\zeta^{\theta,a}_s -)$ exists
     (i.e. is not equal to the cemetery point).
   \item Let $r>s\geq 0$. Conditionally on $\zeta^{\theta,a}$ and
     $W^{\theta,a}_s$, we have $W^{\theta,a}_r$  is equal to
     $W^{\theta,a}_s$ on $[0, m)$, where $m=\min(\zeta^{\theta,a}_s,
     \zeta^{\theta,a}_r)$ and is 
     distributed as $\xi$ under $\rP_{\hat W^{\theta,a}_s}$ on $[m,
     \zeta^{\theta,a}_r)$. 
\end{itemize}
We recall the snake property: a.s. for all $s,s'$, we have 
$W^{\theta,a}_s(r)=W^{\theta,a}_{s'}(r)$, for all $r<\min (\zeta^{\theta,a}_s,
     \zeta^{\theta,a}_{s'})$.

We first check that the family  $(W^{\theta,a}, a>0)$ can be built in a
consistent way.  Let $\bar \cc^c$  be the set  of càdlàg function  defined on
$\R_+$  taking values  in $\W_y$  and where  the life  time  process is
continuous and  lies in $[0,c]$.   Let $a>b>0$. For $\bar  \varphi \in
\bar \cc^a$, with lifetime  process $\varphi$. Recall the function
$C_\varphi$ defined in Section \ref{sec:hp}. We   define  
the projection, $\Pi_{a,b}$, from 
$\bar \cc^a$ to
$\bar \cc^b$ by  $\Pi_{a,b}(\bar \varphi)=\bar \varphi\circ  C_\varphi$.

Recall \reff{eq:def-T} and define $W^{\theta,a,(x)}=(W^{\theta,a}_s,
s\in [0, T^{\theta,a}_x])$. 

\begin{lem}
\label{lem:Proj}
Let $a>b>0$.    The process $\Pi_{a,b} (W^{\theta,a}) $
(resp. $\Pi_{a,b} (W^{\theta,a,(x)}) $)  is distributed as
   $W^{\theta,b}$ (resp. $W^{\theta,b,(x)} $). 
\end{lem}
The  proof is similar  to the  proof of  Lemma \ref{lem:proj},  see also
Remark  \ref{rem:ch=cst}.   Notice  that  by construction  we  have  for
$a>b>c$,   $\Pi_{a,c}=\Pi_{a,b}\circ   \Pi_{b,c}$.   The   compatibility
relation of  Lemma \ref{lem:Proj} implies the existence  of a projective
limit  $\cw^\theta=(\cw^{\theta,a},a\geq 0)$ such  that $\cw^{\theta,a}$
is distributed as $W^{\theta,a}$ and
\begin{equation}
   \label{eq:PPiH}
\Pi_{a,b} (\cw^{\theta,a})=\cw^{\theta,b}. 
\end{equation}
Similarly,      we      can      define     a      projective      limit
$\cw^{\theta,(x)}=(\cw^{\theta,a,(x)},a\geq           0)$           s.t.
$\cw^{\theta,a,(x)}$   is    distributed   as   $W^{\theta,a,(x)}$   and
\reff{eq:PPiH} holds with  $\cw^{\theta,(x)}$ instead of $\cw^{\theta}$.
The  family of  lifetime  processes of  $(\cw^{\theta,a},  a\geq 0)$  is
distributed as  $\ch^{\theta}$. Therefore, we  shall denote $\ch^\theta$
(resp. $\ch^{\theta,a}$)  the lifetime process  of $\cw^{\theta}$ (resp.
$\cw^{\theta,a}$).    This   notation   is   consistent   with   Section
\ref{sec:RK}.   We call  the  process $\cw^\theta$  the Brownian  snake.
From  Remark  \ref{rem:ch=cst},  notice  that for  $\theta\geq  0$,  the
process $\cw^{\theta,a,(x)}$ is independent of $a$ for $a$ large enough.
We shall identify the projective limit $\cw^{\theta,(x)}$ to this common
value.  It correspond to  the usual Brownian snake in \cite{lg:sbprspde}
when  $\theta=0$ (stopped  when  the local  time  at 0  of its  lifetime
reaches $x$).

\subsection{Excursion and special Markov property}

We denote by $\N_y^{\theta,a}$ the excursion measure of $\cw^{\theta,a}$
away  from  the  trivial  path   $y$,  with  lifetime  $0$.   We  assume
$\N_y^{\theta,a}$ is normalized so  that the corresponding local time at
$y$ (as defined in \cite{b:emp} Chap. 3) is the local time at $0$ of the
lifetime  process:  $L^{\theta,a}_0$. Let  $\sigma^{\theta,a}=\inf\{s>0;
\ch^{\theta,a}_s=0\}$.  Under  $\N^{\theta,a}_y$, $\sigma^{\theta,a}$ is
the length of the lifetime excursion.

\begin{lem}
   We have the first moment formula: for $F$ any non-negative measurable
   function defined on the space of  càdlàg $E$-valued
function
\begin{equation}
   \label{eq:moment1N}
\N^{\theta,a}_y\left[\int_0^{\sigma^{\theta,a}} ds\;
  F(\cw_s^{\theta,a})\right]
= \int_0^a dr\; \expp{-4\theta r} \rP_y\left[ F((\xi_t, t\in
  [0,r]))\right]. 
\end{equation}
\end{lem}
This result is known for $\theta\geq 0$ (see \cite{dlg:rtlpsbp}
proposition 1.2.5). 
\begin{proof}
Let $a>0$ and $x>0$ be fixed. Notice it is enough to prove establish
\reff{eq:moment1N} with $W^{\theta,a}$ instead of $ \cw^{\theta,a}$. 
To be concise, we shall omit $a$ and $x$,
so that we write for example $H^\theta$ for
$H^{\theta,a}$ or $T^\theta$ for $T^{\theta,a}_x$. 

Note that $\{s;  H ^{\theta}_s>0, s\in [0,T^{\theta}]\}$ is open,
and consider $(\alpha_i,\beta_i)$, $i\in I$, its connected
component. Let $G$ be any non-negative measurable
   function defined on the space of  càdlàg $\W_y$-valued
function and set $G_i=G( W ^{\theta}_t, t\in (\alpha_i, \beta_i))$. 

Similar arguments used for the proof of  \reff{eq:GirH} relying on
Girsanov theorem implies that 
\[
\E_y\left[\expp{    -\sum_{i\in I} G_i}    \right]
= \E_y\left[\expp{ \theta x - \theta  L^0_a(T^0)  - 2 \theta^2 T^0} \expp{
    -\sum_{i\in I} G_i} \right].
\]
Excursion theory gives that 
\[
\E_y\left[\expp{    -\sum_{i\in I} G_i}    \right]
= \exp\left\{-x \N_y^\theta[1- \expp{-G( W ^\theta)}]\right\}
\]
 and since  $T^0=\sum_{i\in I} \beta_i-\alpha_i$ and $L^0_a(T^0)=\sum_{i\in
   I} L^0_a(\beta_i) - L^0_a(\alpha_i)$, 
\[
\E_y\left[\expp{ \theta x - \theta  L^0_a(T^0)  - 2 \theta^2 T^0} \expp{
    -\sum_{i\in I} G_i} \right]
=
\exp{\left\{\theta x -x \N_y^0[1- \expp{ - \theta
    L^0_a(\sigma^0) -2\theta^2 \sigma^0-G( W ^0)}]\right\}}. 
\]
Thus, we get 
\[
\N_y^\theta[1- \expp{-G( W ^\theta)}]
=
-\theta  + \N_y^0[1- \expp{- \theta
    L^0_a(\sigma^0)-2\theta^2 \sigma^0  -G( W ^0)}].
\]
First moment computation implies that 
\begin{equation}
   \label{eq:NG1}
\N_y^\theta [G( W ^\theta)]
= 
\N_y^0 [ G( W ^0) \expp{ - \theta
    L^0_a(\sigma^0)- 2\theta^2 \sigma^0 }]. 
\end{equation}
Now we specialize to the case $G( W ^\theta)=\int_0^{\sigma^\theta}
F( W ^\theta_s)\; ds$, where $F$ is a  
non-negative measurable
   function defined on the space of  càdlàg $E$-valued
function. 
We have
\[
\N_y^\theta \left[\int_0^{\sigma^\theta} F( W ^\theta_s)\; ds\right]
= 
\N_y^0 \left[ \int_0^{\sigma^0} F( W ^0_s) \expp{ - \theta
    L^0_a(\sigma^0)- 2\theta^2 \sigma^0 }\; ds\right]. 
\]
We  can replace  $\displaystyle  \expp{  -  \theta
  L^0_a(\sigma^0)-  2\theta^2  \sigma^0 }$  in  the  right   side
member  by  its  predictable 
projection  
\[
   \expp{      -   \theta
  L^0_a(s)-    2\theta^2   s }\E\left[\expp{ -\theta
    L^0_a(\sigma^0)- 2\theta^2 \sigma^0 }\mid H^0_0=r\right]_{\mid r=H^0_s}.
\]
Using  time
reversibility (see  Corollary 3.1.6 
of \cite{dlg:rtlpsbp} in a more general case) for the first equality and
predictable projection for the second equality, we get that
\begin{multline*}
    \N_y^0 \left[ \int_0^{\sigma^0} F( W ^0_s) \expp{  - \theta
    L^0_a(s)- 2\theta^2 s } \E\left[\expp{ - \theta
    L^0_a(\sigma^0)- 2\theta^2 \sigma^0}\mid H^0_0=r\right]_{\mid r=H^0_s} \;
ds\right]\\ 
\begin{aligned}
&=
\N_y^0 \left[ \int_0^{\sigma^0} F( W ^0_s) \expp{  - \theta
    (L^0_a(\sigma) -L^0_a(s))- 2\theta^2 (\sigma- s)} \E\left[\expp{ -\theta
    L^0_a(\sigma^0)- 2\theta^2 \sigma^0  }\mid H^0_0=r\right]_{\mid
  r=H^0_s} \; ds\right]\\ 
&=
  \N_y^0 \left[ \int_0^{\sigma^0} F( W ^0_s) \E\left[\expp{ -\theta
    L^0_a(\sigma^0)- 2\theta^2 \sigma^0}\mid H^0_0=r\right]_{\mid
  r=H^0_s}^2 \; ds\right]. 
\end{aligned}
  \end{multline*}
Notice from \reff{eq:sde} that 
\[
\inv{2} L^0_a(\sigma^0)
= H_{\sigma^0} - \beta_{\sigma^0}
- \inv{2} L^0_0(\sigma^0) - H_0
= - \beta_{\sigma^0}- H_0.
\]
Stopping time theorem for exponential martingale implies that 
\[
\E\left[\expp{ -\theta
    L^0_a(\sigma^0)- 2\theta^2  \sigma^0}\mid H^0_0=r\right]
=\E\left[\expp{ -2\theta \beta_{\sigma^0} -
    2\theta H^0_0- 2\theta^2  \sigma^0 }\mid H^0_0=r\right]
=\expp{-2\theta r}. 
\]
We deduce that 
\[
\N_y^\theta \left[\int_0^{\sigma^\theta} F( W ^\theta_s)\; ds\right]
= 
\N_y^0 \left[ \int_0^{\sigma^0} F( W ^0_s) \expp{ - 4 \theta H^0_s }\;
  ds\right].  
\]
The result is then a consequence of the first moment formula for the
Brownian snake (see formula (4.2) in \cite{dlg:rtlpsbp}). 

\end{proof}

We can define the exit local time of an open subset of $D$ of $E$. For $w \in \W$, let $\tau(w)=\inf\{t>0; w(t)\not\in
D\}$. Let  $y\in D$. We assume that $\rP_y(\tau<\infty )>0$. Following
\cite{lg:sbprspde, lg:bssd}, the limit 
\[
L^D_s=\lim_{\varepsilon \downarrow 0} \inv{\varepsilon} \int_0^t
\ind_{\{
  \tau(\cw^{\theta,a}_u)<\zeta_u<\tau(\cw^{\theta,a}_u)+\varepsilon\}}\;
du 
\]
exists for all $s>0$ $\P_y^{\theta}$-a.s. and $\N_y^{\theta,a}$-a.e. and
defines a  continuous non-decreasing additive functional.  We deduce the
first moment formula from \reff{eq:moment1N} and proposition 4.3.2 in
\cite{dlg:rtlpsbp}:
\begin{equation}
   \label{eq:moment1Ntau}
\N^{\theta,a}_y\left[\int_0^{\sigma^{\theta,a}} dL^D_s\;
  F(\cw_s^{\theta,a})\right]
=   \rP_y\left[ \expp{-4\theta \tau} \ind_{\{\tau\leq a \}} F((\xi_t, t\in
  [0,\tau]))\right]. 
\end{equation}

We consider the following hypothesis:

$(\ca)$ For every $y\in D$, $w$ is continuous at $s=\tau$, 
$\rP_y$-a.s. on $\{\tau<\infty \}$. 

\noindent
We recall the description of the excursions of $\cw^{\theta,a}$ out of
$D$. We consider
\begin{equation}
   \label{eq:def-A}
A_t=\int_0^t \ind_{\{ \ch^{\theta,a}_s\leq  \tau(\cw^{\theta,a}_s)\}}\;
ds,
\end{equation}
and $\eta_s=\inf\{t; A_t>s\}$ its right continuous inverse. We define
the càdlàg process $\tilde
\cw^{a}_s=\cw_{\eta_s}^{\theta,a}$. Let $\tilde \cf=(\tilde \cf_t, t\geq
0)$ be 
the filtration generated by $\tilde  \cw^{a}$. Note the snake
property implies that $\{s\in [0, \sigma^{\theta,a}];
\ch^{\theta,a}_s>  \tau(\cw^{\theta,a}_s) \}$ is open and consider
$(\alpha_i, \beta_i)$, $i\in I$, its connected component. The excursions
of $\cw^{\theta,a}$ out of $D$, $\cw^{\{i\}}$, $i\in I$,  are defined by
\[
\cw^{\{i\}}(r)=\cw^{\theta,a}_{(\alpha_i+t)\wedge
  \beta_i}(r+\ch^{\theta,a}_ {\alpha_i}), \quad r\in [0, \ch^{\{i\}}_t=\ch^{\theta,a}_{(\alpha_i+t)\wedge
  \beta_i}].
\]
We denote by $\sigma^{\{i\}}=\alpha_i-\beta_i$ the duration of the
excursion $\cw^{\{i\}}$. Following the proof of theorem 2.4 of
\cite{lg:bssd}, one can  check the next result.

Let $\D$ the space of càdlàg function defined on $\R_+$ taking values  in
$\W$ and let $\delta_{z}$ denote the Dirac mass at point
$z$. 
\begin{prop}
\label{prop:exit-measure}
   The random measure $\int_0^{\sigma^{\theta,a}} dL^D_s \;
   \delta_{(L^D_s, \ch^{\theta,a}_s)}$ is measurable w.r.t. $\tilde
   \cf_\infty $. Let $\phi$ a non-negative measurable function defined on
   $\R_+\times \R_+\times \D$, we have 
\[
\N_y^{\theta,a} \left[\expp{- \sum_{i\in I} \phi(L^D_{\alpha_i},
    \ch^{\theta,a}_{\alpha_i}, \cw^{\{i\}})}\mid \!\tilde \cf_\infty \right]
= \exp\left\{ - \int_0^{\sigma^{\theta,a}}  d L^D_s\;
  \N^{\theta,a-h}_{\hat \cw^{\theta,a}_{s}} [1-\expp{-
    \phi(\ell, h,\cdot)}]_{| \ell= L^D_s, h=\ch^{\theta,a}_{s}}\!\right\}.
\]
\end{prop}

\section{Pruning of the height process}
\label{sec:p}

We present a pruning of the genealogical tree described by the height
process $\ch^\theta$ using a method introduced in \cite{as:psf} in the
case $\theta=0$. The pruning gives a natural way to recover
$Z^{\theta+\gamma}$ from $Z^\theta$ for $\gamma>0$. This gives a dual
procedure to \cite{ad:cbpimcsbpi}, where the authors used immigration
to reconstruct $Z^\theta$ from $Z^{\theta+\gamma}$. This pruning procedure goes
back to \cite{ds:scabs}, where the authors used an intensity of the
killing rate which was dependent of the underlying motion. 

\subsection{Poisson process as  underlying motion}
We keep notations of the previous Section.

Let  $\gamma >0$, and following  \cite{as:rbsd,as:psf} consider  for the
spatial motion (i.e.  the law of $\xi$) the Poisson process distribution
with intensity $4\gamma $.  
We denote by $L^D$ the exit local time out of
$D=\{0\}$.  The additive functional $A$ defined in \reff{eq:def-A} 
can be written in the following way:
\[
A_t=\int_0^t \ind_{\{ \cw^{\theta,a}_s=0\}}\;
ds
\]
Let  $\tilde \ch^a$ be the  lifetime process of $\tilde \cw^a$. 
Notice that $(\tilde \ch^a, a\geq 0)$  is  a consistent  family in  the sense
that $\pi_{a,b}  ( \tilde \ch^a)=\tilde \ch^b$ for  all $a\geq
b\geq 0$. We shall denote by $\tilde \ch$ its projective limit and call
it the pruned height process.

\begin{lem}
\label{lem:A=L}
We have a.s. $\P_0^{\theta,a}$-a.s. and $\N_0^{\theta,a}$-a.e.  for all
$s\geq 0$, $L_s^D=4\gamma  A_s$. 
\end{lem}

\begin{proof}
We shall first prove the result for $\theta=0$. We drop the notation
$a$ and $\theta=0$ in the first part of the proof. We have 
\begin{align*}
   \N_0[(L_\sigma^D -L_t^D-4 \gamma  A_\sigma+4 \gamma  A_t)^2]
&=2\N_0\left[\int_t^\sigma (L_\sigma^D -L_s^
D - 4\gamma  A_\sigma + 4\gamma 
  A_s) d(L_s^D -4\gamma  A_s)\right]\\
&=2\N_0\left[\int_t^\sigma  2 \ch_s\N[L^D_\sigma -4\gamma  A_\sigma]  d(L_s^D -4\gamma  A_s)\right],
\end{align*}
where we  used the previsible  projection of $(L_\sigma^D -L_s^D  - 4\gamma 
A_\sigma  + 4\gamma   A_s)$  and proposition  2.1  in \cite{lg:bssd}  to
compute  it   for  the  second  equality.   Now  \reff{eq:moment1N}  and
\reff{eq:moment1Ntau}                    implies                    that
$\N_0[L_\sigma^D]=4\gamma \N_0[A_\sigma]$.      This      implies     that
$\N_0[(L_\sigma ^D -L_t^D-4 \gamma   A_\sigma+4 \gamma   A_t)^2]=0$  for all
$t\geq 0$. Since $L^D$ and $A$ are continuous and equal to $0$ at $0$,
this implies  that $\N_0$-a.e.  for  all $s\geq 0$,  $L_s^D=4\gamma  A_s$.
Since  $\P_0$-a.s.    $\int_0^\sigma     \ind_{\{\ch_s=0\}}    dL_s^D    =\int_0^\sigma
\ind_{\{\ch_s=0\}} dA_s =  0$, we deduce from excursion  theory that the
result holds also $\P_0$-a.s.

Using   Girsanov  theorem   (see   \reff{eq:NG1}),  since   $L^D=4\gamma A$
holds $\N^{0,a}_0$-a.e,   we    deduce   that   the    equality   also   holds
$\N^{\theta,a}_0$-a.e. (and also $\P_0^{\theta,a}$-a.s.).

\end{proof}

Using Lemma \ref{lem:A=L} notice that $ \N^{\theta,a}_0$-a.e. 
\[
\int_0^{\sigma^{\theta,a}} dL^D_s \;
   \delta_{(L^D_s, \ch^{\theta,a}_s)}
=4 \gamma  \int_0^{\tilde \sigma^{a}} du \;
   \delta_{(u, \tilde \ch^{a}_u)},
\]
where  $\tilde \sigma^{a}=\inf\{s>0; \tilde
\ch^a_s=0\}=A_{\sigma^{\theta,a}}$  is the  length of  the excursion  of $\tilde
\cw^{a}$. Let
us also notice that $\hat \cw^a_{s}=0$ $dA_s$-a.e.

The Poisson process  does not satisfy condition $(\ca )$ with
$D=\{0\}$. However Proposition \ref{prop:exit-measure}  can be  extended
to  this particular  case. (See  \cite{ad:falp} for a  similar formulation in
slightly  different context.   In \cite{ad:falp},  there is  no Brownian
part, and  the Poisson process  is only increasing  at the nodes  of the
height process.) Using the previous remarks, 
Proposition \ref{prop:exit-measure} can be written as follows.

\begin{prop}
\label{prop:MS}
Let $\phi$ a non-negative measurable function defined on
   $\R_+\times \R_+\times \D$, we have 
\[
\N_0^{\theta,a} \left[\expp{- \sum_{i\in I} \phi(A_{\alpha_i},
    \ch^{\theta,a}_{\alpha_i}, \cw^{\{i\}})}\mid \tilde \cf_\infty \right]
= \exp\left\{ - 4\gamma  \int_0^{\tilde \sigma^{a}}  d r\;
  \N^{\theta,a-h}_0 [1-\expp{-
    \phi(r,h,\cdot)}]_{| h=\tilde \ch^{a}_{r}}\right\}.
\]
\end{prop}

\subsection{The main result}

\begin{prop}
\label{prop:pruning}
The   pruned   height   process   $\tilde   \ch$   is   distributed   as
$\ch^{\theta+\gamma }$.
\end{prop}

Recall that  the height  process $\ch^{\theta}$ allows  to code  for the
genealogy of continuous state branching process with branching mechanism
$\psi_\theta$.   In  fact,  using  the Poisson  process  with  intensity
$4\gamma $ as  a spatial motion provides a way  to remove individuals of
continuous  state branching  process  associated to  the height  process
$\ch^{\theta}$ in such a way  as to preserve the genealogical structure.
The  height process  corresponding  to the  remaining  individuals is  a
height    process     associated    to    the     branching    mechanism
$\psi_{\theta+\gamma  }$.  This  Proposition  is  an  extension  to  the
super-critical case of \cite{as:psf}.

\begin{proof}
   Because of the consistency, it is enough to prove that $\tilde \ch^a$
   is distributed as $H^{\theta+\gamma ,a}$. Let $\theta\in  \R$, $a\geq
   0$ be fixed. We shall omit $\theta$ and
   $a$ in what follows and for example write $\cw_s$ for
   $\cw^{\theta,a}_s$. 

   Recall  $\eta$  is  the   right  continuous  inverse  of  $A$,  where
   $\displaystyle A_t=\int_0^t \ind_{\{ \cw^{\theta,a}_s=0\}}\; ds$. 
We shall use a sub-martingale problem, see \cite{sv:dpbc}, to give the
law of $\tilde \ch$.  
Recall that $\ch $ solves \reff{eq:sde}:
\[
\ch_t= \beta_t - 2\theta t + \inv{2} L_0(t) -\inv{2} L_a(t),
\]
where $\beta=(\beta_t, t\geq 0)$ is a Brownian motion. Let $g$ be defined
on $[0,\infty )\times \R$ with compact support with first derivative in
the first variable  and second derivative in the second variable
continuous (in both variables). We shall write $g'(t,x)=\partial_x
g(t,x) $, $g''(t,x)=\partial^2_{xx} g(t,x)$. We shall assume that
$g'(t,0)\geq 0$ and $g'(t,a)\leq 0$ for all $t\geq 0$. 
We define for $t\geq
0$ 
\[
M_t=g(0,0) + \int_0^t g'(A_s,\ch_s) \; d\beta_s+\inv{2}\int_0^t
g'(A_s,0) dL_0(s) -\inv{2}\int_0^t
g'(A_s,a)\ind_{\{\hat \cw_s =0\}} dL_a(s).
\]
Notice that $(M_t, t\geq 0)$ is a sub-martingale with respect to
$\cf=(\cf_t, t\geq 0)$,
the filtration generated by $\ch$. We also have 
\begin{multline}
   \label{eq:MartgHt}
M_t=g(A_t,\ch_t)-\int_0^t \left(\inv{2}
  g''(A_s,\ch_s) - 2\theta g'(A_s,\ch_s)\right)\; ds -\int_0^t
\partial_t g(A_s, \ch_s) dA_s \\
\nonumber +\inv{2}\int_0^t
g'(A_s,a)\ind_{\{\hat \cw_s \neq 0\}} dL_a(s).
\end{multline}
Since  $\eta_t$ is an $\cf$-stopping
time, the  stopping time Theorem implies the process
 $N=(N_t, t\geq 0)$, where $N_t=\E[M_{\eta_t}| \tilde \cf_t]$,  is an
$\tilde \cf$-sub-martingale. 
We set 
\[
\tilde M_t= \int_0^{\eta_t} \left(\inv{2} g''(A_s,\ch_s) - 2\theta
g'(A_s,\ch_s)\right)\ind_{\{\hat \cw _s\neq
0\}}\; ds 
- \frac{1}{2}\int_0^{\theta_t}g'(A_s, a) \ind_{\{\hat\cw_s\neq 0\}}\; d
L_a(s) . 
\]
Recall that a.s.  $A_{\eta_t}=t$ to get 
\begin{align*}
   M_{\eta_t}
&= g(t,\tilde \ch_t) - \int_0^{\eta_t} \left(\inv{2} g''(A_s,\ch_s) - 2\theta
g'(A_s, \ch_s)+  \partial_t g(A_s, \ch_s)
\right)\ind_{\{\hat\cw_s=0\}}\; ds 
   -\tilde M_t \\
&= g(t,\tilde \ch_t) - \int_0^{t} \left(\inv{2} g''(s,\tilde \ch_s) - 2\theta
g'(s,\tilde \ch_s)+  \partial_t g(s, \ch_s) \right)\; ds 
 -\tilde M_t. 
\end{align*}
We have, using notations of Proposition \ref{prop:MS},
\begin{align*}
   \tilde M_t
&=  \int_0^{\eta_t} \left(\inv{2} g''(A_s, \ch_s) - 2\theta
g'(A_s, \ch_s)\right)\ind_{\{\hat \cw_s\neq 0\}}\; ds
-\frac{1}{2}\int_0^{\eta_t} g'(A_s, a)
  \ind_{\{\hat\cw_s\neq 
  0\}}\; d L_a(s)\\ 
&=  \sum_{i\in I} \ind_{\{A_{\alpha_i}\leq
    t\}}\int_{0}^{\sigma^i} \left(\inv{2} 
  g''(A_{\alpha_i}, \ch_s^i+\ch_{\alpha_i} ) - 2\theta
g'(A_{\alpha_i}, \ch_s^i+ \ch_{\alpha_i} )\right)\; ds\\
&\hspace{6cm}
- \frac{1}{2} \sum_{i\in I} \ind_{\{A_{\alpha_i}\leq
    t\}} g'(A_{\alpha_i}, a)  (L_a(\beta_i)
  -L_a(\alpha_i))  . 
\end{align*}
We get
\begin{multline*}
 \E\left[ \tilde M_t| \tilde \cf_\infty \right]\\  
\begin{aligned}
&= 4\gamma  \int_0^t du \;  \N^{a-h}_0\left[ \int_0^\sigma \left(\inv{2} 
  g''(u,\ch_s+h) - 2\theta
g'(u,\ch_s+h)\right)\; ds -\frac{g'(u,a)}{2}L_a(\sigma) \right]_{|h=\tilde
\ch_ u} \\ 
&= 4\gamma  \int_0^t du \left[\int_0^ {a-\tilde \ch_u}  \expp{-4\theta s }
\left(\inv{2}  
  g''(u,s+\tilde \ch_u) - 2\theta
g'(u,s+\tilde \ch_u) \right) \; ds - \frac{g'(u,a)}{2} \expp{-4\theta (a-\tilde
\ch_u)} \right]  \\ 
&= 4\gamma  \int_0^t du  \left[\left[\inv{2} g'(u,s+\tilde \ch_u)\expp{-4\theta
    s }\right]_0^{a-\tilde \ch_u} - \frac{g'(u,a)}{2}\expp{-4\theta
  (a-\tilde \ch_u)}\right]\\
&=- 2\gamma  \int_0^t du \; g'(u,\tilde \ch_u),
\end{aligned}
\end{multline*}
where we used  Proposition \ref{prop:MS} for the first  equality and 
\reff{eq:moment1N} and  \reff{eq:moment1Ntau} for the second. This
implies 
\[
 \E\left[\tilde M_t| \tilde \cf_t \right]=- 2\gamma 
\int_0^t du \; g'(u,\tilde \ch_u), 
\] 
and we deduce that  
\begin{align*}
N_t
&=\E[M_{\eta_t}| \tilde \cf_t]\\
&= g(t,\tilde \ch_t) - \int_0^{t} \left(\inv{2} g''(s,\tilde \ch_s) - 2\theta
g'(s,\tilde \ch_s)+ \partial_t g(s, \tilde \ch_s) \right)\; ds  - \E\left[\tilde M_t| \tilde \cf_t \right]\\
&= g(t,\tilde \ch_t) - \int_0^{t} \left(\inv{2} g''(s,\tilde \ch_s) -
  2(\theta+\gamma ) 
g'(s,\tilde \ch_s)+ \partial_t g(s, \tilde \ch_s) \right)\; ds.
\end{align*}
Notice  that a.s.  $\tilde  \ch_t\in  [0,a]$. Recall  $N$  is a  $\tilde
\cf$-sub-martingale for  any smooth function $g$  such that $g'(t,0)\geq
0$ and $g'(t,a)\leq 0$ for all  $t\geq 0$.  We deduce from uniqueness of
solution to the sub-martingale  problem, see \cite{sv:dpbc} theorem 5.5,
that $\tilde  \ch$ is distributed as  a Brownian motion  in $[0,a]$ with
drift $-2(\theta+\gamma  )$ and  reflected at 0  and $a$. This  and the
consistency property end the
proof.
\end{proof}

\newcommand{\sortnoop}[1]{}

\end{document}